\documentclass[11pt,fleqn,twoside]{article}
\usepackage{amsfonts,amssymb,latexsym}
\makeatletter
\newcommand{\prava}{\footnotesize\it
\begin{flushright}
\begin{minipage}{18cm}
Copyright \copyright 1998 by P.E. Hydon
\end{minipage}
\end{flushright}}

\newcommand{\name}[1]{\begin{flushleft}
                       \LARGE \bf #1
                       \end{flushleft}\vspace{-3mm}}

\newcommand{\Author}[1]{\begin{flushleft}
                       \it #1 \end{flushleft}}

\newcommand{\Adress}[1]{\begin{flushleft}
                       \it #1 \end{flushleft}}

\newcommand{\Date}[1]{\begin{flushleft}
                      \small  \it #1 \end{flushleft}}

\newcommand{\ehkol}{Author \ name}
\newcommand{\ohkol}{Article \ name}
\renewcommand{\@evenhead}{
\hspace*{-3pt}\raisebox{-15pt}[\headheight][0pt]{\vbox{\hbox to \textwidth
{\thepage \hfil \ehkol}\vskip4pt \hrule}}}
\renewcommand{\@oddhead}{
\hspace*{-3pt}\raisebox{-15pt}[\headheight][0pt]{\vbox{\hbox to \textwidth
{\ohkol \hfil \thepage}\vskip4pt\hrule}}}
\renewcommand{\@evenfoot}{}
\renewcommand{\@oddfoot}{}

\newcommand{\be}{\begin{equation}}
\newcommand{\ee}{\end{equation}}
\newcommand{\ba}{\hspace*{-5pt}\begin{array}}
\newcommand{\ea}{\end{array}}

\newcommand{\ds}{\displaystyle}
\makeatother

\begin{document}
\setcounter{page}{405}
\thispagestyle{empty}

\renewcommand{\ehkol}{P.E. Hydon}
\renewcommand{\ohkol}{How to Find Discrete Contact Symmetries}

\begin{flushleft}
\footnotesize \sf
Journal of Nonlinear Mathematical Physics \qquad 1998, V.5, N~4,
\pageref{hydon-fp}--\pageref{hydon-lp}.
\hfill {\sc Article}
\end{flushleft}

\vspace{-5mm}

\renewcommand{\footnoterule}{}
{\renewcommand{\thefootnote}{} \footnote{\prava}}

\name{How to Find Discrete Contact Symmetries}\label{hydon-fp}

\Author{Peter E. HYDON}

\Adress{Department of Mathematics and Statistics, \\
 University of Surrey, Guildford GU2 5XH, UK\\
 E-mail: P.Hydon@surrey.ac.uk}

\Date{Received May 15, 1998; Accepted July 7, 1998}

\begin{abstract}
\noindent
This paper describes a new algorithm for determining all discrete contact
symmetries of any dif\/ferential equation whose Lie contact symmetries are known.
The method is constructive and is easy to use. It is based upon the observation
that the adjoint action of any contact symmetry is an automorphism of the Lie
algebra of generators of Lie contact symmetries. Consequently, all contact
symmetries satisfy
various compatibility conditions. These conditions enable the discrete
symmetries to
be found systematically, with little ef\/fort.
\end{abstract}

\section{Introduction}

\setcounter{equation}{0}

Discrete symmetries of dif\/ferential equations are used in various ways.
They map solutions to (possibly new) solutions. They may be used to create
ef\/f\/icient numerical methods for the computation of solutions to boundary-value
problems. Indeed, there is currently much research into techniques for
constructing numerical methods that preserve various types of symmetry
[2, 4, 10].
Discrete and continuous groups of symmetries determine the nature of
bifurcations in
nonlinear dynamical systems. Equivariant bifurcation theory describes the
ef\/fects of
symmetries, but it may yield misleading results unless {\it all} symmetries of
the dynamical system are known [3, 6].

In general, it is straightforward to f\/ind all one-parameter Lie groups of
symmetries
of a given system, using techniques developed by Sophus Lie more than a
century ago [1, 11, 12, 14]. Yet, until recently, no
simple method for f\/inding all discrete symmetries was known. Ansatz-based
methods can be used to f\/ind discrete symmetries belonging to particular
classes, e.g. [5],
but such methods cannot guarantee that all discrete symmetries have been found.
The main dif\/f\/iculty is that, commonly, the determining equations for discrete
symmetries
form a highly-coupled nonlinear system. Reid and co-workers have developed
a computer algebra package aimed at reducing this system to a dif\/ferential
Gr\"obner basis~[13], but the method is computationally intensive and seems
not to have been widely used.

A new approach to the problem of f\/inding discrete point symmetries has recently
been described by the author [7]. Instead of trying to solve the symmetry
condition directly,
one f\/irst examines the
adjoint action of an arbitrary discrete point symmetry upon the Lie algebra of
Lie point
symmetry generators. This yields a set of necessary conditions which simplify
the problem of
constructing all discrete point symmetries.
Of course, one must know the Lie algebra of Lie point symmetry generators,
which
must be non-trivial. However, this is not a severe limitation, as most
dif\/ferential equations of physical importance have a non-trivial Lie algebra.

Most techniques of Lie symmetry analysis readily generalize to contact
symmetries,
and the aim of the current paper is to show how all discrete contact
symmetries of
a given dif\/ferential equation may be
found systematically, using an extension of the algorithm for determining
discrete
point symmetries. Perhaps surprisingly, many dif\/ferential equations without
non-point Lie contact symmetries
have non-point discrete contact symmetries. It is useful to be able to derive
these (non-obvious)
symmetries systematically because, like all symmetries, they constrain the
behaviour of solutions of
the dif\/ferential equation.

The new method is described f\/irst in the context of ordinary dif\/ferential
equations   (ODEs) of order $n\geq 3$, then adapted
to treat a f\/inite subgroup of the inf\/inite group of contact symmetries of a
given second-order
ODE. Finally, it is shown that the method can be extended to partial
dif\/ferential equations (PDEs).

\section{The algorithm applied to ODEs}
\setcounter{equation}{0}

A dif\/feomorphism
\[
\Gamma:(x,y)\mapsto(\hat x,\hat y)
\]
is a symmetry of the ODE
\be
y^{(n)}=\omega\left(x,y,y',\dots,y^{(n-1)}\right)
\ee
if it maps the set of solutions to itself, i.e. if
\be
\hat y^{(n)}=\omega\left(\hat x,\hat y,\hat y',\dots,\hat y^{(n-1)}\right)
\qquad {\rm when\ (2.1)\ holds}.
\ee
Here the functions $\hat y^{(k)}$ are obtained by prolonging the dif\/feomorphism
$\Gamma$ to derivatives, using
\be
\hat y^{(k)}\equiv\frac{D\hat y^{(k-1)}}{D\hat x},\qquad \left(\hat
y^{(0)}=\hat y\right),
\ee
where
\[
D=\partial_{x}+y'\partial_{y}+y''\partial_{y'}+\cdots.
\]
For point symmetries, $\hat x$ and $\hat y$ are functions of $x$ and $y$.
Contact symmetries are more general than point symmetries, because $\hat x$,
$\hat y$ and $\hat y'$ are functions of $x$, $y$ and~$y'$. One-parameter
Lie groups of contact symmetries are obtained in a similar way to Lie point
symmetries, by linearizing the symmetry condition (2.2). Specif\/ically,
\[
\ba{l}
\hat x = x + \epsilon\xi(x,y,y')+O\left(\epsilon^2\right),\\[2mm]
\hat y = y + \epsilon\eta(x,y,y')+O\left(\epsilon^2\right),\\[2mm]
\hat y^{(k)} = y^{(k)} + \epsilon\eta^{(k)}\left(x,y,y',\dots,y^{(k)}\right)+
O\left(\epsilon^2\right),\qquad k\geq 1.
\ea
\]
The contact conditions require that the $O(\epsilon)$ terms should be
expressible in
terms of the characteristic function $Q(x,y,y')=\eta-y'\xi$, as follows
\be
\xi=-Q_{y'},\qquad\eta=Q-y'Q_{y'},\qquad\eta^{(k)}=D^kQ-y^{(k+1)}Q_{y'},\quad
k\geq 1.
\ee
In particular,
\[
\eta^{(1)}(x,y,y')=Q_x+y'Q_y.
\]
Lie point symmetries of ODEs of order $n\geq 2$ are Lie contact symmetries
whose characteristic
function is linear in $y'$.

The set of all inf\/initesimal generators of Lie contact symmetries of a given
ODE of order
$n\geq 3$ forms a f\/inite-dimensional Lie algebra, ${\cal L}$, which can
generally be determined
systematically [14]. Given such an ODE, suppose that ${\cal L}$ has a basis
\be
X_i=\xi_i(x,y,y')\partial_x+\eta_i(x,y,y')\partial_y+\eta^{(1)}_i(x,y,y')
\partial_{y'},
\qquad i=1,\dots,N,
\ee
where $\left(\xi_i,\eta_i,\eta_i^{(1)}\right)$
are obtained from the characteristic
function $Q_i(x,y,y')$
using (2.4), and $N={\rm dim}\,({\cal L})$.
The structure constants, $c_{ij}^k$, for the basis (2.5)
are determined by
\be
[X_i,X_j]=c_{ij}^k X_k.
\ee
(Summation is implied when an index occurs twice, one raised and once lowered.)
The one-parameter Lie group of contact symmetries corresponding to a
particular $X_i$
is obtained by exponentiation. We use the notation
\[
\Gamma_i(\epsilon):(x,y,y')\mapsto\left(e^{\epsilon X_i}x,\,e^{\epsilon
X_i}y,\,e^{\epsilon X_i}y'\right).
\]
Suppose that
\be
\Gamma:(x,y,y')\mapsto \left(\hat x(x,y,y'),\, \hat y(x,y,y'),\, \hat
y'(x,y,y')\right)
\ee
 is a contact symmetry of the given ODE.
 Then the contact transformation obtained  by the
adjoint action of $\Gamma$ upon $\Gamma_i(\epsilon)$,
\be\hat \Gamma_i (\epsilon)=\Gamma \Gamma_i(\epsilon)\Gamma^{-1}, \ee is
also a contact symmetry, for each $\epsilon$ in some neighbourhood of zero.
Therefore, for each~$i$, there is a (local) one-parameter Lie group of contact
symmetries
\be
\hat\Gamma_i(\epsilon):(\hat x,\hat y,\hat y')\mapsto
\left(e^{\epsilon \hat X_i}\hat x,\,e^{\epsilon \hat X_i}\hat y,\,
e^{\epsilon \hat X_i} \hat y'\right),
\ee
whose inf\/initesimal generator is
\be
\hat X_i=\Gamma X_i\Gamma^{-1}=\xi_i(\hat x,\hat y,\hat y')\partial_{\hat
x}+\eta_i(\hat x,\hat y,
\hat y')\partial_{\hat y}+\eta_i^{(1)}(\hat x,\hat y,\hat y')\partial_{\hat
y'}.
\ee
Consequently
\[
\hat X_i\in{\cal L},\qquad i=1,\dots, N.
\]
The generators $\{\hat X_i\}_{i=1}^N$ are simply the basis generators
$\{X_i\}_{i=1}^N$ with $(x,y,y')$ replaced by $(\hat x,\hat y,\hat y')$.
Therefore
the set $\{\hat X_i\}_{i=1}^N$ is a basis for ${\cal L}$, and so each $X_i$
can
be written as a linear combination of the $\hat X_j$'s. Also, the mapping
$X_i\mapsto\hat X_i$ is an
automorphism of ${\cal L}$ which preserves all structure
constants, i.e.
\be
[\hat X_i,\hat X_j]=c_{ij}^k\hat X_k\qquad {\rm when\ (2.6)\ holds.}
\ee
These results generalize to partial dif\/ferential equations, and are summarized
as follows.

\medskip

\noindent
{\bf Lemma 1.} {\it
Every contact symmetry $\Gamma$ of an ordinary differential equation of order
$n\geq 3$ induces an
automorphism of the Lie algebra, ${\cal L}$,
of generators of one-parameter
local Lie groups of contact symmetries of the differential equation.
For each such $\Gamma$, there exists a constant non-singular matrix
$\left(b_i^l\right)$ such that
\be
X_i=b_i^l\hat X_l.
\ee
This automorphism preserves all structure constants.}

\medskip

Lemma 1 yields the following PDEs for the unknown functions
$\hat x(x,y,y')$ and $\hat y(x,y,y')$:
\be
X_i\hat x=b_i^l\hat X_l\hat x=b_i^l\xi_l(\hat x,\hat y,\hat y'),\qquad
i=1,\dots,N,
\ee
\be
X_i\hat y=b_i^l\hat X_l\hat y=b_i^l\eta_l(\hat x,\hat y,\hat y'),\qquad
i=1,\dots,N.
\ee
This set of $2N$ f\/irst-order PDEs, together with the contact condition
\[
\hat y'(x,y,y')=\frac{\ds \frac{d\hat y}{dx}}{\ds \frac{d\hat x}{dx}},
\]
provides necessary, but not suf\/f\/icient, conditions for $\Gamma$ to be a
contact symmetry.
The contact condition yields the following pair of PDEs, because $\hat y'$ is
independent of $y''$:
\be
\hat y_x+y'\hat y_y=(\hat x_x+y'\hat x_y)\hat y',\qquad \hat y_{y'}=\hat
x_{y'}\hat y'.
\ee
N.B. The lemma gives a further $N$ PDEs
\[
X_i\hat y'=b_i^l\hat X_l\hat y'=b_i^l\eta_l^{(1)}(\hat x,\hat y,\hat y'),\qquad
i=1,\dots,N,
\]
but these add nothing new, for they are a consequence of (2.13), (2.14) and the
contact condition.

To f\/ind all discrete contact symmetries, proceed as follows. First solve the
system of
PDEs (2.13), (2.14), to obtain $(\hat x,\hat y)$
in terms of $x,y,y',b_i^l$ and some unknown constants (or functions) of
integration.
If $N$ is suf\/f\/iciently large, it may be possible to solve this system
algebraically
(by eliminating the derivative terms); otherwise, the method of
characteristics should
be used.
Incorporate the contact condition (2.15); this generally reduces the number of
candidate
solutions $(\hat x,\hat y)$. Finally, use the symmetry condition (2.2) to
determine which
of these solutions are symmetries. The continuous symmetries may be factored
out at a convenient
point in the calculation. The remaining discrete symmetries are inequivalent
under any
continuous symmetry, and form a discrete (but not necessarily f\/inite) group.

If ${\cal L}$ is non-abelian, some of the structure constants are non-zero,
enabling the matrix $B=\left(b_i^l\right)$ to be simplif\/ied {\it
before} any of the 
above calculations are done. Substituting~(2.12) into (2.6), and
taking (2.11) 
into account,
we obtain the following constraints on the components of $B$:
\be
c_{lm}^nb_i^lb_j^m=c_{ij}^kb_k^n.
\ee
It is suf\/f\/icient to
restrict attention to equations (2.16) with $i<j$, because the structure
constants are antisymmetric in the two lower indices.
Moreover, at least some of the continuous symmetries can be factored out
using their adjoint action upon the generators in ${\cal L}$, (see [11]),
\be
{\rm Ad}(\exp(\epsilon_j X_j))X_i=X_i-\epsilon_j[X_j,X_i]+{\epsilon_j^2\over
2!}[X_j,[X_j,X_i]]-\cdots
=a_i^p(\epsilon_j,j)X_p.
\ee
Let $A(j)$ denote the matrix whose components are
$a_i^p(\epsilon_j, j)$, as def\/ined by (2.17).
The system (2.12) is equivalent, under the group generated by $X_j$, to
\[
X_i=\tilde b_i^l\hat X_l,
\]
where $\tilde b_i^l$ are the components of
\be
\tilde B=A(j)B.
\ee
The mapping $B\mapsto\tilde B$ does not af\/fect (2.16),
so we will drop tildes as soon as each equivalence transformation has been
made.
Each generator $X_j$ is used in turn to simplify the form of $B$.

To illustrate this procedure, consider the two-dimensional non-abelian Lie
algebra ${\mathfrak a}(1)$,
with a basis $\{X_1$, $X_2\}$ such that
\be
[X_1,X_2]=X_1.
\ee
The only non-zero structure constants are
\[
c_{12}^1=-c_{21}^1=1.
\]
Therefore (2.16) gives
\[
b_1^1b_2^2-b_1^2b_2^1=b_1^1,\qquad 0=b_1^2,
\]
and hence
\[
B=\left[\matrix{b_1^1&0\cr b_2^1&1\cr}\right],\qquad\ b_1^1\neq 0.
\]
The matrices representing the adjoint action of the continuous group on the
Lie algebra are
\[
A(1)=\left[\matrix{1&0\cr -\epsilon_1&1\cr}\right],\qquad
A(2)=\left[\matrix{e^{\epsilon_2}&0\cr 0&1\cr}\right].
\]
Applying (2.18), f\/irst with $j=1$,
$\ds \epsilon_1={b_2^1\over b_1^1}$, then with
$j=2$, $\epsilon_2=-\ln|b_1^1|$, we obtain
\be
B=\left[\matrix{\alpha&0\cr 0&1\cr}\right],\qquad{\rm where}\quad
\alpha\in\{-1,1\}.
\ee
The reduced form of the matrix $B$ is specif\/ic to this particular
Lie algebra, and is independent of the ODE whose Lie point
symmetries are generated by the algebra.

Simplif\/ied matrices for other non-abelian Lie algebras can be found by the
same technique.
However, if ${\cal L}$ is abelian, the entries of $B$ cannot be determined
{\it a priori}.

\section{Examples}
\setcounter{equation}{0}
The third-order ODE
\be
y'''=\frac{y''^2}{x}-\frac{y''}{y'}
\ee
has a two-dimensional Lie algebra of generators of Lie contact symmetries.
These Lie symmetries are
actually
point symmetries, and ${\cal L}$ is isomorphic to ${\mathfrak a}(1)$.
The basis of~${\cal L}$,
\[
X_1=\partial_y,\qquad X_2=\frac{x}{2}\partial_x+y\partial_y
+\frac{y'}{2}\partial_{y'},
\]
has the commutation relations (2.19), and therefore $B$ is given by (2.20).
The system of PDEs (2.13), (2.14)
amounts to
\[
\left[\matrix{X_1\hat x&X_1\hat y\cr X_2\hat x&X_2\hat
y\cr}\right]=\left[\matrix{\alpha&0\cr 0&1\cr}\right]\left[
\begin{array}{cc}
0&1\\[2mm]
\ds \frac{\hat x}{2}&\hat y \end{array}\right]
=\left[\begin{array}{cc}
0&\alpha\\[2mm] \ds \frac{\hat x}{2}&\hat y\end{array}
\right].
\]
The general solution of this system is
\[
\hat x=xp(t),\qquad\hat y=\alpha y+x^2q(t),\qquad{\rm where}\quad
t=\frac{y'}{x}.
\]
The contact condition gives
\[
\hat y'=xr(t),
\]
where
\[
\alpha t+2q-t\dot q=(p-t\dot p)r\qquad{\rm and}\qquad \dot q =\dot p r.
\]
(Here a dot over a function denotes its derivative with respect to $t$.)
Re-arranging these conditions, we obtain
\be
q=\frac{1}{2}\left( pr-\alpha t\right),
\ee
and the contact condition is satisf\/ied if and only if
\be
p\dot r-\dot p r=\alpha.
\ee
It is convenient to work in terms of $p$ and $r$, because the ODE is invariant
under translations in $y$,
and so $y$ does not occur in the symmetry condition. The prolongation to
second and third derivatives is
\[
\hat y''=\frac{r+(y''-t)\dot r}{p+(y''-t)\dot p},
\]
\[
\hat y'''=\frac{(p\dot r-\dot p r)xy'''+(y''-t)^2(p\ddot r-\ddot p
r)+(y''-t)^3(\dot p\ddot r-\ddot p\dot r)}{x(p+(y''-t)\dot p)^3}.
\]
These expressions are substituted into the symmetry condition (2.2), and
powers of $y''$ are equated
to yield an over-determined system of nonlinear ODEs for $p$ and $r$. These
are easily solved with the
aid of the contact condition (3.3). There are two sets of solutions. Either
\be
(p,r)=(c,ct),\qquad c^2=\alpha,
\ee
or
\be
(p,r)=(ct,c),\qquad c^2=-\alpha.
\ee
Re-writing these solutions in terms of the original variables, we obtain eight
inequivalent
discrete symmetries that form a group isomorphic to ${\mathbb Z}_4\times
{\mathbb Z}_2$. The
group generators are
\[
\Gamma_1:(x,y,y')\mapsto(ix,-y,iy'),
\]
\[
\Gamma_2:(x,y,y')\mapsto(y',xy'-y,x).
\]
The subgroup generated by $\Gamma_1$ consists of four discrete point
symmetries, which can be found
without having to consider contact symmetries [7]. The four inequivalent
non-point contact
symmetries are obtained from the point symmetries by composition with
$\Gamma_2$, which
is the prolonged Legendre Transform.

Generally speaking, the larger the dimension of ${\cal L}$, the easier it is
to f\/ind
the discrete symmetries. Nevertheless, the contact condition makes it possible
to solve
the governing equations, even when $N=1$. Suppose that ${\cal L}$ is
one-dimensional and
the ODE is written in canonical coordinates, so that
the continuous symmetries
are generated by
\[
X=\partial_y.
\]
Then Lemma 1 gives $X=b\hat X$, where $b\neq 0$, and therefore
\[
\hat x=f(x,y'),\qquad\hat y=by+g(x,y').
\]
The contact condition yields
\[
\hat y'=h(x,y'),
\]
where
\be
g_x=f_xh-by', \qquad g_{y'}=f_{y'}h.
\ee
Equations (3.6) are compatible if and only if
\be
f_xh_{y'}-f_{y'}h_x=b.
\ee
They can be integrated, once $f$ and $h$ are known, to determine $g$ up to an
arbitrary
constant, which may be set at any convenient value to factor out equivalence
under the
one-parameter group generated by $X$.

Consider the general third order ODE admitting the group generated by $X$:
\[
y'''=\omega(x,y',y'').
\]
Substituting
\[
\hat x=f(x,y'),\qquad\hat y'=h(x,y')
\]
into the symmetry condition, and then equating powers of $y''$, yields an
over-determined coupled
system of nonlinear PDEs. This system is precisely as intractible as the
problem of using the
symmetry condition alone to f\/ind all point symmetries of
\[
y''=\omega(x,y,y').
\]
However, with the aid of the contact condition (3.7), the problem simplif\/ies
considerably.

To illustrate this, consider the ODE
\be
y'''=y''^3\sin\left(\frac{x}{y'}\right),
\ee
whose only Lie contact symmetries are those generated by $X=\partial_y$.
The symmetry condition gives the over-determined (but complicated) system
\[
f_xh_{xx}-f_{xx}h_x=h_x^3\sin\left(\frac{f}{h}\right),
\]
\[
f_{y'}h_{xx}+2f_{x}h_{xy'}-f_{xx}h_{y'}-2f_{xy'}h_{x}=3h_x^2h_{y'}\sin
\left(\frac{f}{h}\right),
\]
\[
f_{x}h_{y'y'}+2f_{y'}h_{xy'}-f_{y'y'}h_x-2f_{xy'}h_{y'}=3h_xh_{y'}^2\sin
\left(\frac{f}{h}
\right),
\]
\[
f_{y'}h_{y'y'}-f_{y'y'}h_{y'}+(f_{x}h_{y'}-f_{y'}h_{x})\sin\left(\frac{x}{y'}
\right)=h_{y'}^3\sin\left(\frac{f}{h}\right).
\]
This system can be greatly simplif\/ied by using (3.7) and its dif\/ferential
consequences, which
reduces the f\/irst three equations to
\[
f_{xx}=f_{xy'}=h_x=h_{y'y'}=0.
\]
Combining this result with (3.7) and the remaining symmetry condition gives
\[
f=\alpha(x+2n\pi y'),\qquad h=y',\qquad \alpha\in\{-1,1\},\quad n\in{\mathbb Z}.
\]
After solving (3.6) for $g(x,y')$ and setting $g(0,0)=0$ to factor out the
continuous symmetries,
we obtain the following result. The inequivalent discrete contact symmetries
of (3.8) form a
countably-inf\/inite group, which is generated by
\[
\Gamma_1:(x,y,y')\mapsto(-x,-y,y'),
\]
\[
\Gamma_2:(x,y,y')\mapsto(x+2\pi y',y+\pi y'^2, y').
\]

\section{Discrete uniform contact symmetries}
\setcounter{equation}{0}

For ODEs of order $n\geq 3$, the Lie contact symmetries can generally be found
systematically, and
${\cal L}$ is f\/inite-dimensional. Second-order ODEs have an
inf\/inite-dimensional Lie algebra of
contact symmetry generators, but they cannot all be found unless the general
solution of the ODE is
known. However, some contact symmetries may be found with the aid of a
suitable ansatz for $Q$.
For example, the restriction $Q_{y'y'}=0$ enables all Lie point symmetries to
be found systematically.

Other restrictions on $Q$ are possible. For example, the set of all {\it
uniform} contact symmetries
of a given ODE of order $n\geq 2$ is a f\/inite-dimensional Lie group [8].
Uniform contact symmetries
are of the form
\be
\hat x=\Phi(x,y'), \qquad\hat y=ky+\Theta(x,y'),
\qquad\hat y'=\Psi(x,y'),\qquad
k\in{\mathbb R}\backslash\{0\},
\ee
where the contact condition requires that
\be
ky'+\Theta_x=\Phi_x\Psi,\qquad \Theta_{y'}=\Phi_{y'}\Psi.
\ee
Many of the most commonly-occurring contact symmetries are uniform, including
all contact symmetries
of the ODEs that were used as examples in the previous section.
The generators of uniform Lie contact symmetries have characteristic functions
of the form
\[
Q=c_1y+\phi(x,y'),\qquad c_1\in {\mathbb R}.
\]
The Lie algebra of these generators can be found systematically, for a given
ODE of order $n\geq 2$,
by equating powers of $y$ in the symmetry condition. The constructions leading
to Lemma 1 can be
repeated, restricting attention to uniform contact symmetries, to obtain the
following.

\medskip

\noindent
{\bf Lemma 2.}  {\it
Every uniform contact symmetry $\Gamma$ of an ordinary differential equation
of order $n\geq 2$ induces an
automorphism of the Lie algebra, ${\cal L}$, of generators of one-parameter
local Lie groups of uniform contact symmetries of the differential equation.
For each such $\Gamma$, there exists a constant non-singular matrix
$\left(b_i^l\right)$ such that
\be
X_i=b_i^l\hat X_l.
\ee
This automorphism preserves all structure constants.}

\medskip

This lemma can be used to f\/ind all discrete uniform contact symmetries in
exactly the same way as
Lemma 1 is used. To illustrate this, consider the ODE
\be
y''=\frac{y'^2}{3xy'-4y},
\ee
which has a four-dimensional Lie algebra of uniform contact symmetry
generators. The structure
constants are simplest in the basis
\[
X_1=2y'\partial_x+y'^2\partial_y,
\]
\[
X_2=\frac{3}{4}x\partial_x+\frac{1}{2}y\partial_y-\frac{1}{4}y'\partial_{y'},
\]
\[
X_3=2(y')^{-3}\partial_x+3(y')^{-2}\partial_y,
\]
\[
X_4=\frac{1}{4}x\partial_x+\frac{1}{2}y\partial_y+\frac{1}{4}y'\partial_{y'}.
\]
The only non-zero structure constants are
\[
c^1_{12}=-c^1_{21}=1,\qquad c^3_{34}=-c^3_{43}=1,
\]
and so the Lie algebra is isomorphic to
${\mathfrak a}(1)\oplus{\mathfrak a}(1)$.
After simplifying $B$ as
far as possible, using the relations (2.16) and the adjoint action of the
continuous group, we obtain
two possibilities. Either
\be
B=\left[\matrix{\alpha &0&0&0\cr 0&1&0&0\cr 0&0&\beta&0\cr
0&0&0&1\cr}\right],\qquad\ \alpha,\beta\in\{-1,1\},
\ee
or
\be
B=\left[\matrix{0&0&\alpha &0\cr 0&0&0&1\cr \beta&0&0&0\cr
0&1&0&0\cr}\right],\qquad\ \alpha,\beta\in\{-1,1\}.
\ee
The inequivalent discrete uniform contact symmetries are calculated in the
same way as previously. Lemma 2
is used, together with (4.1) and the contact condition (4.2), to obtain the
most general form possible
for a uniform contact symmetry of (4.4). The symmetry condition for the ODE is
then used to determine
which of the possible solutions actually are symmetries. We f\/ind that the
group of inequivalent
discrete uniform symmetries of~(4.4) is isomorphic to
${\mathbb Z}_2\times{\mathbb Z}_2\times{\mathbb Z}_2$, and is
generated by
\[
\Gamma_1:(x,y,y')\mapsto(-x,y,-y'),
\]
\[
\Gamma_2:(x,y,y')\mapsto(x,-y,-y'),
\]
\[
\Gamma_3:(x,y,y')\mapsto\left(xy'^2,2xy'-y,(y')^{-1}\right).
\]
The discrete point symmetries generated by $\Gamma_1$ and $\Gamma_2$ can also
be found
directly from the Lie algebra of point symmetries, which is ${\rm
Span}(X_2,X_4)$. These
inequivalent discrete symmetries are derived from the matrix $B$ in (4.5),
with $\alpha=\beta$;
they map each of the ${\mathfrak a}(1)$ Lie subalgebras to itself. However
$\Gamma_3$, which is
derived from (4.6), interchanges these two subalgebras.
Actually, (4.4) is merely one representative of a whole class of second-order
ODEs that have
inequivalent uniform discrete symmetries $\Gamma_i$,
$i=1,2,3;$ these ODEs are of the form
\[
y''= \frac{my'^2}{y-(m+1)xy'},\qquad m\in{\mathbb R}\backslash\{0\}.
\]

\section{Contact symmetries of a nonlinear PDE}
\setcounter{equation}{0}

The method described in section 2 generalizes to PDEs without dif\/f\/iculty; the
corresponding
algorithm for point symmetries will be discussed elsewhere [9]. The basic
steps are shown
here, using the potential hyperbolic heat equation
\be
u_{tt}+u_t=\frac{u_{xx}}{u_x}
\ee
as an example. The Lie algebra of contact symmetry generators is
f\/ive-dimensional, with
a basis
\[
X_1=\partial_t,\qquad X_2=\partial_x,\qquad X_3=\partial_u,\qquad
X_4=e^{-t}\partial_u-e^{-t}\partial_{u_t},
\]
\[
X_5=-x\partial_x+u\partial_u+u_t\partial_{u_t}+2u_x\partial_{u_x}.
\]
These Lie contact symmetries are actually point symmetries, but there is still
the possibility
that some discrete contact symmetries may be non-point symmetries, {\it cf.}
(3.1). Using
(2.16) and the adjoint action of the continuous group to simplify $B$, we
obtain two possibilities. Either
\be
B={\rm diag}\,\{1,b,\alpha,\beta,1\},\qquad\alpha,\beta\in\{-1,1\},\quad
b\in{\mathbb R}\backslash\{0\},
\ee
or
\be
B=\left[\matrix{1&0&0&0&0\cr 0&0&b&0&0\cr 0&\alpha&0&0&0\cr 0&0&0&\beta&0\cr
2&0&0&0&-1\cr}\right],\qquad \alpha,\beta\in\{-1,1\},\quad b\in{\mathbb
R}\backslash\{0\}.
\ee
If $B$ is of the form (5.2), the equations analogous to (2.13), (2.14) are
\[
\left[\matrix{X_1\hat t&X_1\hat x&X_1\hat u\cr X_2\hat t&X_2\hat x&X_2\hat
u\cr
X_3\hat t&X_3\hat x&X_3\hat u\cr X_4\hat t&X_4\hat x&X_4\hat u\cr X_5\hat
t&X_5\hat x&X_5\hat u\cr}\right] =B\ \left[\matrix{1&0&0\cr 0&1&0\cr 0&0&1\cr
0&0&e^{-\hat t}\cr 0&-\hat x&\hat u\cr}
\right],
\]
whose general solution is
\[
\hat t=t+c_1,\qquad\hat x=bx+c_2|u_x|^{-\frac{1}{2}},\qquad \hat u=\alpha
u+(\alpha-\beta
e^{-c_1})u_t+c_3|u_x|^{\frac{1}{2}},\quad c_i\in{\mathbb R}.
\]
The contact condition and the symmetry condition reduce these further, to
\[
(\hat t, \hat x,\hat u)=(t, \alpha x,\alpha u),\qquad\alpha\in\{-1,1\}.
\]
The remaining inequivalent discrete contact symmetries are obtained from (5.3)
in a similar
way. To summarize the results: the inequivalent discrete contact symmetries of
(5.1)
form a group isomorphic to ${\mathbb Z}_2\times{\mathbb Z}_2$, which is generated by
\[
\Gamma_1:(t,x,u,u_t,u_x)\mapsto(t,-x,-u,-u_t,u_x),
\]
\[
\Gamma_2:(t,x,u,u_t,u_x)\mapsto\left(t+\ln|u_x|,\, u+u_t,\,
x+u_t(u_x)^{-1},-u_t(u_x)^{-1},(u_x)^{-1}\right).
\]
The symmetry $\Gamma_2$ was known previously [15], but it was not known that
$\Gamma_2$ and $\Gamma_1\Gamma_2$
are the only non-point contact symmetries (up to equivalence). The method
outlined in the current paper
enables the user to completely classify all contact symmetries of a given
dif\/ferential equation with a known non-trivial Lie algebra.

 \label{hydon-lp}

\end{document}